\documentclass[12pt]{article}
\usepackage{amsfonts}
\usepackage{amssymb}

\input{tcilatex}
\begin{document}

\begin{center}
\textbf{Variants on the Berz sublinearity theorem}

\bigskip

by

\bigskip

\textbf{N. H. Bingham and A. J. Ostaszewski}\\[0pt]
\bigskip

\textit{To Roy O. Davies on his 90}$^{\text{\textit{th}}}$ \textit{birthday.}

\bigskip
\end{center}

\noindent \textbf{Abstract. }We consider variants on the classical Berz
sublinearity theorem, using only DC, the Axiom of Dependent Choices, rather
than AC, the Axiom of Choice which Berz used. We consider thinned versions,
in which conditions are imposed on only part of the domain of the function
-- results of quantifier-weakening type. There are connections with
classical results on subadditivity. We close with a discussion of the
extensive related literature.

\bigskip

\noindent \textbf{Keywords: }Dependent Choices, Homomorphism extension,
Subadditivity, Sublinearity, Steinhaus-Weil property, Thinning, Quantifier
weakening.\textbf{\ \newline
}


\bigskip

\textbf{1. Introduction: sublinearity.}

We are concerned here with two questions. The first is to prove, as directly
as possible, a linearity result via an appropriate group-homomorphism
analogue of the classical Hahn-Banach Extension Theorem $HBE$ [Ban] -- see
[Bus] for a survey. Much of the $HBE$ literature most naturally elects as
its context \textit{real} Riesz spaces (ordered linear spaces equipped with
semigroup action, see \S 4.8), where some naive analogues can fail -- see
[BusR]. These do not cover our test-case of the additive reals $\mathbb{R}$,
with focus on the fact (e.g. [BinO3]) that for $\mathbb{A\subseteq R}$ a
dense subgroup, if $f:\mathbb{A}\rightarrow \mathbb{R}$ is additive (i.e. a
partial homomorphism) and locally bounded, then it is linear: $f(x):=ca$ for
some $c\in \mathbb{R}$ and all $a\in \mathbb{A}$. Can this result be deduced
by starting with some natural, continuous, subadditive majorant $S:\mathbb{R}%
\rightarrow \mathbb{R}$ (so that, equivalently, $\bar{S}|\mathbb{A\leq }%
f\leq S|\mathbb{A}$ for $\bar{S}(.)=-S(-.),$ which is super-additive) and
then invoking an (interpolating) additive extension $F$ majorized by $S?$
For then $F,$ automatically being continuous, is linear, because its
restriction to the rationals $F|\mathbb{Q}$ is so (as in Th. 1 below).
Assuming additionally positive-homogeneity, $HBE$ yields an $F$, but this
strategy relies very heavily on powerful selection axioms (formally a
weakend version of the Prime Ideal Theorem, itself a weakening of the Axiom
of Choice, AC, see \S 4.7). The alternative is to apply either semigroup
results in [Kau], [Fuc], [Kra], or the recent group-theoretic result in
[Bad3], but all these again rely on AC (see \S 4.7 again). We give an answer
in Theorem 3 that relies on the much weaker axiom of Dependent Choices, DC
(see \S 4.7 once more). We stress that, throughout the paper, \textit{all
our results need only DC}.

The group analogue (for $\mathbb{R}$) of \textit{sublinearity} used by [Ber]
(cf. [Kau]) requires subadditivity as in Banach's result [Ban, \S 2.2 Th.
1], but restricts Banach's positive-homogeneity condition to just $\mathbb{N}
$-homogeneity:%
\[
S(nx)=nS(x)\qquad (x\in \mathbb{R}\text{, }n=0,1,2,...) 
\]%
(with the universal quantifier $\forall $ on $x$ and $n$ understood here, as
is usual in mathematical logic). From here onwards we take this to be our
definition of sublinearity. This is of course equivalent to
`positive-rational-homogeneity'. Berz proves and uses a Hahn-Banach theorem
in the context of $\mathbb{R}$ as a vector space over $\mathbb{Q}$ (for
which see also [Kuc2, \S 10.1]) to show that if $S:\mathbb{R}\rightarrow 
\mathbb{R}$ is measurable and sublinear, then $S|\mathbb{R}_{+}$ and $S|%
\mathbb{R}_{-}$ are both linear; for generalizations to Baire (i.e. having
the \textit{Baire} property) and universally measurable functions in
contexts including Banach spaces, again using only DC, see [BinO2]. Berz's
motivation was questions of \textit{normability} in topological spaces
[Ber]. The key result here is Kolmogorov's theorem [Kol]: normability is
equivalent to the origin having a bounded convex neighbourhood ([Rud, Th.
1.39 and p. 400]).

Our second, linked, question asks whether the universal quantifier $(x\in 
\mathbb{R})$ above can be \textit{weakened} to range over an additive
subgroup. Since $S=\mathbf{1}_{\mathbb{R}\backslash \mathbb{Q}}$ is
subadditive and $\mathbb{N}$-homogeneous on $\mathbb{Q}$, but not linear on $%
\mathbb{R}_{\pm },$ the \textit{quantifier weakening} must be accompanied by
an appropriate side-condition. We give in Theorem 5 a necessary and
sufficient condition (referring also to a thinned-out domain), by extending
the standard asymptotic analysis -- as in [HilP] (see Theorem HP below) --
of the ratio $S(t)/t$ near $0$ and at infinity; this, indeed, permits
thinning-out the universal quantifier of $\mathbb{N}$-homogeneity to a dense
additive subgroup $\mathbb{A}$.

\bigskip

We come at these questions here employing ideas on quantifier weakening
previously applied in [BinO3] to additivity issues in classical regular
variation, and in [BinO2] to Jensen-style convexity in Banach spaces. We
borrow from [BinO3] two key tools: Theorem 0 below on continuity (exploiting
an idea of Goldie), and Theorem 0$^{+}$ on linear (upper) bounding
(exploiting early use by Kingman of the Baire Category Theorem -- see
[BinO1]), the latter delayed till \S 3, when we have the preparatory results
needed.

\bigskip

\noindent \textbf{Theorem 0. }\textit{For subadditive} $S:\mathbb{R}%
\rightarrow \mathbb{R}\cup \{-\infty ,+\infty \}$ \textit{with }$%
S(0+)=S(0)=0:$\textit{\ }$S$\textit{\ is continuous at }$0$ \textit{iff }$%
S(z_{n})\rightarrow 0,$\textit{\ for some sequence }$z_{n}\uparrow 0,$%
\textit{\ and then }$S$ \textit{is continuous everywhere, if finite-valued.}

\bigskip

In \S 2 below we discuss subadditivity, sublinearity and theorems of Berz
type, proving Theorems 1-3, Th. BM (for Baire/measurable) and Th. HP (for
Hille-Phillips). The work of Hille and Phillips is a major ingredient in the
Kingman subadditive ergodic theeorem (\S 4.9) of probability theory. In \S 3
we give stronger versions of Berz's theorem by thinning the domain of
definition, under appropriate side-conditions, results of
quantifier-weakening type (Theorems 4 and 5). We close in \S 4 with a
discussion of the extensive background literature.

\bigskip

\textbf{2. Subadditivity, sublinearity and theorems of Berz type.}

We begin with a sharpened form of the Berz theorem, with a proof that seems
new. Here and below we write $B_{\delta }(x):=(x-\delta ,x+\delta )$ for the
open $\delta $-ball around $x.$

\bigskip

\noindent \textbf{Theorem 1 }(cf. [BinO3]\textbf{).} \textit{For} $S:\mathbb{%
R}\rightarrow \mathbb{R}$ \textit{a sublinear function (i.e. subadditive,
with }$S(nx)=nS(x)$ \textit{for }$x\in \mathbb{R}$\textit{\ and} $%
n=0,1,2,...)$\textit{, if }$S$ \textit{is locally bounded, then both }$S|%
\mathbb{R}_{+}$ \textit{and} $S|\mathbb{R}_{-}$ \textit{are linear.}

\bigskip

\noindent \textbf{Proof. }For $M$ a bound on $S$ in $B_{\delta }(0)=(-\delta
,\delta ),$%
\[
|S(x)|=|S(kx)/k|\leq M/k, 
\]%
for $k\in \mathbb{N}$, $x\in B_{\delta /k}(0);$ so $S\ $is continuous at $0,$
and so everywhere. But $S(q)=qS(1)$ for rational $q>0,$ so by continuity $%
S(x)=xS(1)$ for $x\in \mathbb{R}_{+};$ likewise $S(x)=|x|S(-1)$ for $x\in 
\mathbb{R}_{-}.$ $\square $

\bigskip

\noindent \textbf{Remark. }The argument can be repeated for $S:\mathbb{R}%
\rightarrow X$ with $X$ a normed vector space; then $S(x)=||x||S(u_{x})$ for 
$u_{x}$ the unit vector on the ray: $\{\lambda x:\lambda \geq 0\}.$ Here $%
|S(u_{x})|\leq M/\delta $ for all $x\neq 0.$

\bigskip

This gives as a corollary

\bigskip

\noindent \textbf{Theorem BM (}[Ber], [BinO2], cf. [BinO3]\textbf{). }%
\textit{For} $S:\mathbb{R}\rightarrow \mathbb{R}$ \textit{a sublinear
function, if }$S$\textit{\ is Baire/measurable, then both }$S|\mathbb{R}_{+}$
\textit{and} $S|\mathbb{R}_{-}$ \textit{are linear.}

\bigskip

\noindent \textbf{Proof. }For $S$ Baire/measurable, $S$ is bounded above on
a non-negligible set and so, being subadditive, is bounded above on some
interval (by the Steinhaus-Weil Theorem, [Oxt], [BinO2]), and so, being
subadditive, is locally bounded. $\square $

\bigskip

The following is a slightly sharper form of results in [BinO3] with a
simpler proof (the subgroup here is initially arbitrary). This extension
theorem may be interpreted in Hahn-Banach style as involving a subadditive
function $S$ which, relative to a subgroup $\mathbb{A}$, majorizes an
additive function $G$ that happens to agree with the restriction $S|\mathbb{A%
}$.

\bigskip

\noindent \textbf{Theorem 2.} \textit{If} $S:\mathbb{R}\rightarrow \mathbb{R}
$ \textit{is a subadditive locally bounded function and} $\mathbb{A}$ 
\textit{any non-trivial additive subgroup such that }$S|\mathbb{A}$ \textit{%
is additive, then }$S|\mathbb{A}$ \textit{is linear.}

\textit{In particular, for }$\mathbb{A}$\textit{\ dense, any additive
function }$G$ \textit{on }$\mathbb{A}$\textit{\ has at most one continuous
subadditive extension }$S:\mathbb{R}\rightarrow \mathbb{R}$\textit{.}

\bigskip

We will need the following Theorem; for completeness, we show how to make
the simple modification needed to the result given in [HilP, Th. 7.6.1]. (We
replace their additional blanket condition of measurability of $S$ by local
boundedness, and give more of the details, as they are needed later.)

\bigskip

\noindent \textbf{Theorem HP.} \textit{For} $S:\mathbb{R}\rightarrow \mathbb{%
R}$\textit{\ a locally bounded subadditive function }%
\[
\beta =\beta _{S}:=\inf\nolimits_{t>a}\frac{S(t)}{t}=\lim\nolimits_{t%
\rightarrow \infty }\frac{S(t)}{t}<\infty \qquad (a>0), 
\]%
\textit{so }$\beta $\textit{\ does not depend on the choice of }$a>0.$ 
\textit{In particular, }%
\[
\beta _{S}:=\inf\nolimits_{t>0}\frac{S(t)}{t}\in \mathbb{R}. 
\]

\noindent \textbf{Proof. }Following [HilP, Th.7.5.1], for $a>0$ and $ma\leq
t<(m+1)a$ with $m=2,3,...,$ we note two inequalities, valid according as $%
S(a)\geq 0$ or $S(a)<0:$%
\begin{equation}
\frac{S(t)}{t}\leq \frac{mS(a)+S(t-ma)}{t}\leq \frac{S(a)}{a}+\frac{K}{a},%
\text{ if }S(a)\geq 0,\text{ }\frac{S(a)}{2a}+\frac{K}{a},\text{ if }S(a)<0,
\tag{$\dag $}
\end{equation}%
for $K:=\sup |S([0,2a])|$; indeed%
\[
\frac{1}{2a}\leq \frac{1}{a}-\frac{1}{t}\leq \frac{m}{t}\leq \frac{1}{a}. 
\]%
Also $S(t)/t$ itself is bounded on $[a,2a],$ so $\beta =\beta (a)<\infty $
is well-defined.

Suppose first that $\beta >-\infty ,$ and let $\varepsilon >0.$ As $%
\inf\nolimits_{t>0}S(t)/t<\beta +\varepsilon ,$ choose and fix $b\geq a$
with $S(b)/b\leq (\beta +\varepsilon ).$ For any $t\geq 2b$, let $n=n(t)\in 
\mathbb{N}\geq 0$ satisfy $(n+1)b\leq t<(n+2)b;$ then $b<t-nb<2b$ and%
\[
1-\frac{2b}{t}<\frac{nb}{t}\leq 1-\frac{b}{t}. 
\]%
This time with $K=\sup |S([0,2b])|$ a bound on $S$ as above, since $%
S(t)=S(nb+t-nb)\leq S(nb)+S(t-nb),$ 
\[
\beta \leq \frac{S(t)}{t}\leq \frac{nb}{t}\frac{S(b)}{b}+\frac{S(t-nb)}{t}%
\leq \frac{nb}{t}(\beta +\varepsilon )+\frac{K}{t}\leq (\beta +\varepsilon
)+\varepsilon , 
\]%
for $t>\max \{2b,K/\varepsilon \}.$ So $\lim\nolimits_{t\rightarrow \infty
}S(t)/t=\beta .$

The case $\beta =-\infty $ would be similar albeit simpler. In fact it does
not arise. Indeed, writing $T(t)=S(-t),$ which is subadditive and locally
bounded, yields%
\[
\beta _{S}+\beta _{T}=\lim\nolimits_{t\rightarrow \infty }\left[ \frac{S(t)}{%
t}+\frac{T(t)}{t}\right] \geq 0, 
\]%
as $0\leq S(0)\leq S(t)+T(t).$ So $\beta _{S}>-\infty ,$ since $\beta
_{S}\geq -\beta _{T}>-\infty .$ $\square $

\bigskip

\noindent \textbf{Remark.} In fact%
\[
-\beta _{T}=-\lim\nolimits_{t\rightarrow \infty }\frac{T(t)}{t}%
=-\inf\nolimits_{t>0}\frac{T(t)}{t}=\sup\nolimits_{t>0}\frac{S(-t)}{(-t)}%
=\sup\nolimits_{z<0}\frac{S(z)}{z}=\lim\nolimits_{z\rightarrow -\infty }%
\frac{S(z)}{z}. 
\]

\bigskip

\noindent \textbf{Proof of Theorem 2. }Put $G:=$ $S|\mathbb{A}$, and let $%
\beta _{S}$ denote the unique $\beta $ of Th. HP. For any $a\in \mathbb{A}$ $%
\cap (0,\infty ),$ we have, by Th. HP, that\textbf{\ }%
\[
\beta _{S}=\lim\nolimits_{n\rightarrow \infty }\frac{S(na)}{na}=\frac{G(a)}{a%
}. 
\]%
Now for all $a\in \mathbb{A}$, as $G(-a)=-G(a)$ and $G(0)=0$ (by
additivity), $G(a)=\beta _{S}a$. In particular, for $S$ continuous and $%
\mathbb{A}$ dense, $S(t)=\beta _{S}t$ for all $t\in \mathbb{R}$. $\square $

\bigskip

The next extension theorem employs majorization and minorization on a
subspace. The assumption of subgroup divisibility -- $a/k\in \mathbb{A}$ $%
(a\in \mathbb{A}$ , $k\in \mathbb{N)}$ -- is innocuous (as any subgroup may
be extended to a divisible one without change of cardinality). Below we
write $\mathbb{R}_{+}:=[0,\infty ),$ $\mathbb{R}_{-}:=(-\infty ,0],$ and $%
\mathbb{A}_{\pm }:=\mathbb{A\cap R}_{\pm }.$

\bigskip

\noindent \textbf{Theorem 3}. \textit{For }$\mathbb{A}$\textit{\ a dense,
divisible subgroup of }$\mathbb{R}$\textit{\ and a locally bounded sublinear
(in particular, additive) }$S:\mathbb{A}\rightarrow \mathbb{R}$\textit{,}%
\begin{eqnarray*}
S_{\mathbb{A}}^{+}(x) &:&=\lim\nolimits_{\delta \rightarrow 0}\sup
\{S(t):t\in B_{\delta }(x)\cap \mathbb{A}\}, \\
S_{\mathbb{A}}^{-}(x) &:&=\lim\nolimits_{\delta \rightarrow 0}\inf
\{S(t):t\in B_{\delta }(x)\cap \mathbb{A}\},
\end{eqnarray*}%
\textit{define locally bounded, subadditive, indeed sublinear, functions }$%
S_{\mathbb{A}}^{\pm }:\mathbb{R}\rightarrow \mathbb{R}$\textit{\ with}%
\[
S_{\mathbb{A}}^{-}(x)\leq S_{\mathbb{A}}^{+}(x)\qquad (x\in \mathbb{R})\quad 
\text{and}\quad S_{\mathbb{A}}^{-}(a)\leq S(a)\leq S_{\mathbb{A}%
}^{+}(a)\qquad (a\in \mathbb{A}). 
\]%
\textit{Hence}

\noindent (i)\textit{\ }$S_{\mathbb{A}}^{+}|\mathbb{R}_{\pm }=S_{\mathbb{A}%
}^{-}|\mathbb{R}_{\pm };$

\noindent (ii)\textit{\ }$S_{\mathbb{A}}^{\pm }|\mathbb{R}_{\pm }$, $S|%
\mathbb{A}_{\pm }$ \textit{are linear;}

\noindent (iii) $S_{\mathbb{A}}^{\pm }|\mathbb{A}=S$\textit{.}

\noindent (iv)\textit{\ In particular, for }$S$\textit{\ additive, }$S_{%
\mathbb{A}}^{+}=S_{\mathbb{A}}^{-}$ \textit{and is linear, as is also }$S=S_{%
\mathbb{A}}^{+}|\mathbb{A}$\textit{.}

\bigskip

\noindent \textbf{Proof. }To lighten the notation, we write $S^{\pm }$ for $%
S_{\mathbb{A}}^{\pm }.$ Local boundedness of $S^{\pm }$ follows immediately
from local boundedness of $S.$ Subadditivity is routine, and follows much as
in [HilP, \S 7.8]. As regards sublinearity of $S^{\pm }$, note that if $%
a_{n}\rightarrow x$ for $a_{n}\in \mathbb{A}$ with $\lim\nolimits_{n%
\rightarrow \infty }S(a_{n})=S^{\pm }(x),$ then, as $ka_{n}\in \mathbb{A}$
for $k\in \mathbb{N}$, by sublinearity of $S$%
\begin{eqnarray*}
kS^{+}(x) &=&\lim\nolimits_{n\rightarrow \infty
}kS(a_{n})=\lim\nolimits_{n\rightarrow \infty }S(ka_{n})\leq S^{+}(kx), \\
kS^{-}(x) &=&\lim\nolimits_{n\rightarrow \infty
}kS(a_{n})=\lim\nolimits_{n\rightarrow \infty }S(ka_{n})\geq S^{-}(kx).
\end{eqnarray*}%
Similarly, for $k\in \mathbb{N}$, if $a_{n}\rightarrow kx$ for $a_{n}\in 
\mathbb{A}$ with $\lim\nolimits_{n\rightarrow \infty }S(a_{n})=S^{\pm }(kx),$
then $a_{n}/k\rightarrow x$ with $a_{n}/k\in \mathbb{A}$, and so again by
sublinearity of $S,$ 
\begin{eqnarray*}
S^{+}(kx)/k &=&\lim\nolimits_{n\rightarrow \infty }S(a_{n}/k)\leq S^{+}(x),
\\
S^{-}(kx)/k &=&\lim\nolimits_{n\rightarrow \infty }S(a_{n}/k)\geq S^{-}(x).
\end{eqnarray*}%
So $kS^{\pm }(x)=S^{\pm }(kx).$ By Theorem 1 the four functions $S_{\mathbb{A%
}}^{\pm }|\mathbb{R}_{\pm }$ are linear, and so by dominance the two
functions $S|\mathbb{A}_{\pm }$ are continuous at $0$ and so continuous
everywhere (as in Th. 1). So if $a_{n}\rightarrow a$ with $a,a_{n}\in 
\mathbb{A}$, then $\lim\nolimits_{n\rightarrow \infty }S(a_{n})=S(a)=S^{\pm
}(a),$ proving (iii). So $S^{\pm }|\mathbb{A}_{+}=S|\mathbb{A}_{+};$ this
implies $S|\mathbb{A}_{+}\ $is linear and also that $S^{+}|\mathbb{R}%
_{+}=S^{-}|\mathbb{R}_{+}$, since $\mathbb{A}$ is dense, proving (i) and
(ii) on $\mathbb{R}_{+}$ and $\mathbb{A}_{+}$; similarly on $\mathbb{R}_{-}$
and $\mathbb{A}_{-}$. For additive $S$ this means that $S^{+}=S^{-}$ is
linear, as is $S,$ proving (iv). $\square $

\bigskip

\noindent \textbf{Remarks. }Actually $S_{\mathbb{A}}^{+}$ ($S_{\mathbb{A}%
}^{-}$) is upper (lower) semicontinuous, hence Baire. Of course $S^{\pm
}(kx)\leq kS^{\pm }(x),$ so (in view of the first displayed inequality
above, etc.) a less symmetric proof would have fewer steps. Inducing
functions (such as $S^{\pm }$ from $S,$ above) is a method followed
variously, e.g. in [Kuc1, Th. 1], [BinO2, Th. 5].

\bigskip

\textbf{3. Thinning: a stronger Berz theorem.}

We now give a stronger version of the Berz theorem by weakening a condition
of Heiberg-Seneta type by thinning, as in [BinO3], and requiring the
homogeneity assumption to hold on only a dense additive subgroup $\mathbb{A}$
of $\mathbb{R}$; all in all, with rather less than sublinearity, we improve
on Theorem BM. This comes at the price of assuming more about $S.$ To
motivate the next definition, note that for locally bounded subadditive $S,$
the inequalities ($\dag $) of \S 2 imply that for any $a>0$%
\[
\gamma (a):=\sup\nolimits_{t>a}\frac{S(t)}{t}<\infty , 
\]%
as $|S(t)/t|$ is bounded on $[a,2a].$ As $\gamma (a)$ is decreasing for $%
a>0, $ we have%
\[
-\infty <\beta _{S}\leq \lim \sup\nolimits_{t\downarrow 0}\frac{S(t)}{t}\leq
\infty . 
\]%
We note that, with $T$ as in Theorem HP, $\alpha
_{S}:=\sup\nolimits_{z<0}S(z)/z=-\beta _{T}$ is finite (by the remark
above), so the definition below fills the gap for $\sup%
\nolimits_{t>0}S(t)/t, $ by asking apparently a little less.

\bigskip

\noindent \textbf{Definition.} Say that $S:\mathbb{R}\rightarrow \mathbb{R}$
satisfies the \textit{strong Heiberg-Seneta} ($SHS$) condition if%
\begin{equation}
\gamma =\gamma _{S}^{+}:=\lim \sup\nolimits_{t\downarrow 0}\frac{S(t)}{t}%
<\infty .  \tag{$SHS$}
\end{equation}%
See \S 4.4 for the origin of this term. For $S$ subadditive, we will see in
Proposition 2 that this implies its dual:%
\[
-\infty <\gamma _{S}^{-}:=\lim \inf\nolimits_{t\uparrow 0}\frac{S(t)}{t}\leq
\gamma _{S}^{+}. 
\]

Proposition 1, to which we now turn, associates to each subadditive function 
$S$ a sublinear function $S^{\ast }$ dominating $S,$ here and below to be
called the (upper) \textit{sublinear envelope} of $S.$ (Albeit
multiplicatively, [Kau] studies the lower envelope dominated by $S,$ using
instead $S(nx)/n$ -- also noted in [Fuc], cf. [Bad3] -- an approach followed
in [GajK] employing the decreasing sequence $S(2^{n}x)/2^{n}$.) However,
some assumption on $S$ is needed to ensure that $S^{\ast }$ is
finite-valued: recall that the subadditive function $S=\mathbf{1}_{\mathbb{R}%
\backslash \mathbb{Q}}$ is $\mathbb{N}$-homogeneous on $\mathbb{Q}$, yet $%
S^{\ast }=(+\infty )\cdot \mathbf{1}_{\mathbb{R}\backslash \mathbb{Q}}.$

\bigskip

\noindent \textbf{Proposition 1.} \textit{For }$S:\mathbb{R}\rightarrow 
\mathbb{R}$\textit{\ locally bounded and subadditive with }$S(0)=0$\textit{,
the function defined by}%
\[
S^{\ast }(x):=\lim \sup\nolimits_{n\rightarrow \infty }nS(x/n)\qquad (x\in 
\mathbb{R}) 
\]%
\textit{is subadditive and sublinear and dominates }$S$\textit{. If further }%
$S$\textit{\ satisfies }$(SHS),$\textit{\ then for }$t\geq 0.$%
\[
\beta _{S}t\leq S(t)\leq S^{\ast }(t)\leq \gamma _{S}^{+}t. 
\]%
\textit{In particular, }$S(0+)=S^{\ast }(0+)=0;$\textit{\ furthermore, }$%
\gamma _{S}^{+}>-\infty $ \textit{and}%
\[
\sup\nolimits_{t>0}\frac{S(t)}{t}\leq \gamma _{S}^{+}. 
\]

\bigskip

\noindent \textbf{Proof. }By subadditivity of $S,$ for any $n\in \mathbb{N}$%
\[
S(x)=S(n.x/n)\leq nS(x/n)\leq S^{\ast }(x). 
\]%
Evidently $S^{\ast }$ is subadditive (cf. [HilP, 7.2.2, 7.2.3]). Moreover,
as $S(0)=0$, $S^{\ast }$ is $\mathbb{Q}_{+}$-homogeneous, since for fixed $%
k\in \mathbb{N}$%
\begin{eqnarray*}
kS^{\ast }(x) &=&\lim \sup\nolimits_{m\rightarrow \infty }k\cdot mS(kx/km) \\
&\leq &\lim \sup\nolimits_{n\rightarrow \infty }nS(kx/n)\qquad \text{(via
specialization: }n=km\text{)} \\
&=&S^{\ast }(kx)\leq kS^{\ast }(x),
\end{eqnarray*}%
the latter by subadditivity, so that%
\[
S^{\ast }(kx)=kS^{\ast }(x). 
\]

Suppose now that $(SHS)$ holds. Let $\varepsilon >0.$ Then there is $\delta
>0$ with 
\[
S(x)/x\leq \gamma _{S}^{+}+\varepsilon \qquad (0<x<\delta ). 
\]%
Fix $t>0.$ Then for integer $n>t/\delta $%
\[
\frac{S(t/n)}{t/n}\leq \gamma _{S}^{+}+\varepsilon :\qquad nS(t/n)\leq
(\gamma _{S}^{+}+\varepsilon )t, 
\]%
and so taking limsup as $n\rightarrow \infty $%
\[
S^{\ast }(t)\leq (\gamma _{S}^{+}+\varepsilon )t, 
\]%
for $t\geq 0$, as $S^{\ast }(0)=0$ (since $S(0)=0).$ Taking limits as $%
\varepsilon \downarrow 0$ yields%
\[
S^{\ast }(t)\leq \gamma _{S}^{+}t. 
\]%
Furthermore, for $t\geq 0$%
\[
\beta _{S}t\leq S(t)\leq S^{\ast }(t).\qquad \square 
\]

In view of the linear bounding of $S^{\ast }$ (and hence of $S)$ just proved
from $SHS,$ we proceed to a \textit{weaker property} of $S$ in which the
domain of the limsup operation is \textit{thinned-out}. This will
nevertheless also yield linear bounding of $S$ (from above), hence
finiteness of $\gamma _{S}^{+}$, and in turn the bounding of $S^{\ast }.$ We
need a definition and a theorem from [BinO3].

\bigskip

\noindent \textbf{Definition} [BinO3,2,5]. Say that $\Sigma \subseteq 
\mathbb{R}$ is \textit{locally Steinhaus-Weil (SW)}, or has the \textit{SW
property locally}, if for $x,y\in \Sigma $ and, for all $\delta >0$
sufficiently small, the sets%
\[
\Sigma _{z}^{\delta }:=\Sigma \cap B_{\delta }(z), 
\]%
for $z=x,y,$ have the \textit{interior-point property,} that $\Sigma
_{x}^{\delta }\pm \Sigma _{y}^{\delta }$ has $x\pm y$ in its interior. (Here 
$B_{\delta }(x)$ is the open ball about $x$ of radius $\delta .)$ See
[BinO4, Appendix] for conditions under which this property is implied by the
interior-point property of the sets $\Sigma _{x}^{\delta }-\Sigma
_{x}^{\delta }$ (cf. [BarFN]); for a rich list of examples, see \S 4.5. An
obvious example is an open set $\Sigma .$

We now cite from [BinO3] the following result.

\bigskip

\noindent \textbf{Theorem 0}$^{+}$\textbf{.} \textit{Let} $\Sigma \subseteq
\lbrack 0,\infty )$ \textit{be locally SW accumulating at }$0$\textit{.
Suppose }$S:\mathbb{R}\rightarrow \mathbb{R}$\textit{\ is subadditive with }$%
S(0)=0$ \textit{and:}

\noindent $S|\Sigma $ \textit{is linearly bounded above by }$G(x):=cx$%
\textit{, i.e.} $S(\sigma )\leq c\sigma $ \textit{for some }$c$\textit{\ and
all }$\sigma \in \Sigma ,$ \textit{so that in particular,}%
\[
\lim \sup\nolimits_{\sigma \downarrow 0,\text{ }\sigma \in \Sigma }S(\sigma
)\leq 0. 
\]%
\textit{\ Then }$S(x)\leq cx$ \textit{for all }$x>0,$\textit{\ so} 
\[
\lim \sup_{x\downarrow 0}S(x)\leq 0, 
\]%
\textit{and so }$S(0+)=0.$

\textit{In particular, if furthermore there exists a sequence }$%
\{z_{n}\}_{n\in N}$\textit{\ with }$z_{n}\uparrow 0$\textit{\ and }$%
S(z_{n})\rightarrow 0,$\textit{\ then }$S$\textit{\ is continuous at }$0$%
\textit{\ and so everywhere.}

\bigskip

\noindent \textbf{Definition.} Say that $S:\mathbb{R}\rightarrow \mathbb{R}$
satisfies the \textit{weak Heiberg-Seneta} ($WHS$) condition if for some $%
\Sigma \subseteq (0,\infty ),$ a locally SW set accumulating at $0,$%
\[
\gamma _{S}^{\Sigma }:=\lim \sup\nolimits_{t\downarrow 0,t\in \Sigma }\frac{%
S(t)}{t}<\infty . 
\]

\noindent \textbf{Corollary.} \textit{For }$S:\mathbb{R}\rightarrow \mathbb{R%
}$\textit{\ locally bounded and subadditive with }$S(0)=0$\textit{, if }$S\ $%
\textit{satisfies }$WHS,$\textit{\ then }$S$\textit{\ is linearly bounded by 
}$\gamma _{S}^{\Sigma }t$\textit{\ for }$t\geq 0,$ \textit{and so satisfies }%
$SHS$\textit{\ with }$\gamma _{S}^{+}\leq \gamma _{S}^{\Sigma }.$

\bigskip

\noindent \textbf{Proof. }Write $\gamma =\gamma _{S}^{\Sigma }.$ Let $%
\varepsilon >0.$ Then there is $\delta >0$ with 
\[
S(t)\leq (\gamma +\varepsilon )t\qquad (t\in \Sigma \cap (0,\delta )). 
\]%
So $S(t)\leq (\gamma +\varepsilon )t$ for all $t>0,$ by Th. 0$^{+}$ applied
to $c=\gamma +\varepsilon $. Taking limits as $\varepsilon \downarrow 0$
yields $S(t)\leq \gamma t$ for all $t>0$ and so%
\[
\gamma _{S}^{+}:=\lim \sup\nolimits_{t\downarrow 0,t\in \Sigma }\frac{S(t)}{t%
}\leq \gamma _{S}^{\Sigma }<\infty . 
\]%
So $S$ satisfies the $SHS$. $\square $

\bigskip

We now derive in Theorem 5 below a form of Berz's Theorem, in which the weak
Heiberg-Seneta condition on $S$ permits a thinned-out assumption of
homogeneity; the argument is based on the following result, a corollary of
Theorem 1 and Prop 1.

\bigskip

\noindent \textbf{Theorem 4. }\textit{For }$S:\mathbb{R}\rightarrow \mathbb{R%
}$\textit{\ locally bounded and subadditive with }$S(0)=0$\textit{, if }$S$%
\textit{\ satisfies }$WHS$\textit{, and }$S^{\ast }(t_{n})\rightarrow 0$ 
\textit{for some sequence} $t_{n}\uparrow 0,$\textit{\ then }$S$ \textit{and
its sublinear envelope }$S^{\ast }$\textit{\ are continuous, and further, by
sublinearity, both }$S^{\ast }|\mathbb{R}_{+}$ \textit{and} $S^{\ast }|%
\mathbb{R}_{-}$ \textit{are linear.}

\bigskip

\noindent \textbf{Proof. }By the Corollary we may assume that $SHS\ $holds.
By Prop. 1, $S^{\ast }(0+)=0,$ so $S^{\ast }$ is continuous by Theorem 0,
and now the linearity conclusion follows by Theorem 1, as $S^{\ast }$ is
sublinear (and locally bounded at $0,$ so everywhere -- cf. [BinO3, Prop.
5]). In fact, directly so, since, by homogeneity and continuity, $S^{\ast
}(x)=xS^{\ast }(1)$ and $S^{\ast }(-x)=xS^{\ast }(-1)$ for $x>0,$ as $%
\mathbb{Q}$ is dense and $S^{\ast }(\pm q)=qS^{\ast }(\pm 1)$ for $q\in 
\mathbb{Q}_{+}.$

Now for $x\geq 0$, by subadditivity $-S(x)\leq S(-x)$ (as $S(0)=0)$ and so%
\[
-S^{\ast }(x)\leq -S(x)\leq S(-x)\leq S^{\ast }(-x)=xS^{\ast }(-1). 
\]%
So $S(0-)=0,$ as $S^{\ast }$ is continuous; so $S$ is continuous. $\square $

\bigskip

\noindent \textbf{Proposition 2.} \textit{In the setting of Theorem 3, for} $%
t>0$%
\[
S^{\ast }(-t)/(-t)=\gamma _{S}^{-}\leq \alpha _{S}\leq \beta _{S}\leq \gamma
_{S}^{+}=S^{\ast }(t)/t. 
\]%
\textit{In particular,}%
\[
-\infty <\gamma _{S}^{-}\leq \gamma _{S}^{+}<\infty . 
\]

\bigskip

\noindent \textbf{Proof. }As before we may assume that $SHS$ holds. Write $%
\gamma ^{\pm }:=S^{\ast }(\pm t)/(\pm t)$ for $t>0.$ From Prop. 1 $\beta
_{S}\leq \gamma _{S}^{+};$ as $S\leq S^{\ast },$ for $t>0,$ $\gamma
^{+}=S^{\ast }(t)/t\geq S(t)/t,$ so $\gamma _{S}^{+}\leq \gamma ^{+}.$ For
the reverse inequality, take any $\varepsilon >0$ and choose $\delta >0$
such that $S(t)/t\leq \gamma _{S}^{+}+\varepsilon $ for all $0<t<\delta .$
As $S^{\ast }(1)=\gamma ^{+},$ there exists $m>1/\delta $ with $\gamma
^{+}-\varepsilon \leq mS(1/m).$ Taking $t=1/m<\delta $ yields%
\[
\gamma ^{+}-\varepsilon \leq \gamma _{S}+\varepsilon :\qquad \gamma ^{+}\leq
\gamma _{S}+2\varepsilon . 
\]%
Taking limits as $\varepsilon \downarrow 0$ yields $\gamma ^{+}\leq \gamma
_{S}^{+}.$ Combining, $\gamma ^{+}=\gamma _{S}^{+}.$

Now recall from above that $\alpha _{S}\leq \beta _{S}$ (via $\alpha
_{S}=-\beta _{T}).$ We obtain $\gamma ^{-}\leq \alpha _{S}$ from 
\[
\lim\nolimits_{t\rightarrow -\infty }\frac{S^{\ast }(t)}{t}\leq
\lim\nolimits_{t\rightarrow -\infty }\frac{S(t)}{t}=\alpha _{S}, 
\]%
again as $S\leq S^{\ast }$ ($t$ being negative here).

Put $T^{\ast }(t):=S^{\ast }(-t),=\gamma ^{-}(-t)$ for $t>0.$ Then%
\begin{equation}
\gamma ^{-}=\lim \inf\nolimits_{z\uparrow 0}\frac{S^{\ast }(z)}{z}; 
\tag{$\ast $}
\end{equation}%
indeed 
\begin{eqnarray*}
-\gamma ^{-} &=&\lim\nolimits_{t\rightarrow \infty }\frac{T^{\ast }(t)}{t}%
=\lim \sup\nolimits_{t\downarrow 0}\frac{T^{\ast }(t)}{t}=\lim
\sup\nolimits_{t\downarrow 0}\frac{S^{\ast }(-t)}{t} \\
&=&-\lim \inf\nolimits_{t\downarrow 0}\frac{S^{\ast }(-t)}{-t}=-\lim
\inf\nolimits_{z\uparrow 0}\frac{S^{\ast }(z)}{z}.
\end{eqnarray*}%
By $(\ast )$ and the definition of $\gamma _{S}^{-},$%
\[
\gamma _{S}^{-}=\lim \inf\nolimits_{z\uparrow 0}\frac{S(z)}{z}\geq \lim
\inf\nolimits_{z\uparrow 0}\frac{S^{\ast }(z)}{z}=\gamma ^{-}, 
\]%
as $S\leq S^{\ast }$ ($z$ here being negative). So

\[
\gamma _{S}^{-}\geq \gamma ^{-}. 
\]

We now show that $\gamma _{S}^{-}\leq \gamma ^{-}.$ This runs analogously to
the plus version. Let $\varepsilon >0.$ Choose $\delta >0$ with $\gamma
_{S}^{-}-\varepsilon \leq S(t)/t$ for $t\in (-\delta ,0).$ As $-\gamma
^{-}=S^{\ast }(-1)$, pick $m$ with $m>1/\delta $ and $-\gamma
^{-}-\varepsilon \leq mS(-1/m).$ Then taking $t=-1/m$ gives%
\[
\gamma _{S}^{-}-\varepsilon \leq -mS(-1/m)\leq \gamma ^{-}+\varepsilon
:\qquad \gamma _{S}^{-}\leq \gamma ^{-}+2\varepsilon . 
\]%
Taking limits as $\varepsilon \downarrow 0$ yields $\gamma _{S}^{-}\leq
\gamma ^{-}.$ Combining, $\gamma _{S}^{-}=\gamma ^{-}$. $\square $

\bigskip

\noindent \textbf{Remark. }The burden of proof falls on showing that $\gamma
_{S}^{\pm }=\gamma ^{\pm };$ of course, for $t>0,$ $S(t)+S(-t)=\gamma
^{+}t-\gamma ^{-}t\geq 0$ yields directly that $\gamma ^{-}\leq \gamma ^{+}.$

\bigskip

\noindent \textbf{Theorem 5 (Quantifier-weakened Berz Theorem).} \textit{For 
}$S:\mathbb{R}\rightarrow \mathbb{R}$\textit{\ locally bounded and sublinear
(in particular for} $S$ \textit{Baire/measurable and sublinear) with }$%
S(0)=0 $\textit{, if}\newline
\noindent (i)\textit{\ }$S$ \textit{satisfies }$WHS$\textit{\ and,}\newline
\noindent (ii) $\mathbb{A}$ \textit{is a dense additive subgroup of }$%
\mathbb{R}$ \textit{with }$S|\mathbb{A}$ $\mathbb{N}$-\textit{homogeneous}%
\newline
\noindent \textit{-- then both }$S|\mathbb{R}_{+}$ \textit{and} $S|\mathbb{R}%
_{-}$ \textit{are linear: for }$t\geq 0:$%
\[
S(t)=\beta _{S}t,\text{ and }S(-t)=-\alpha _{S}t. 
\]%
\textit{In particular, for }$S$\textit{\ additive}%
\[
S(t)=\beta _{S}t\qquad (t\in \mathbb{R}). 
\]

\bigskip

\noindent \textbf{Proof. }By the Corollary we may assume that $SHS\ $holds.
Consider any $a\in \mathbb{A}.$ Then $S^{\ast }(a)=\lim
\sup\nolimits_{n\rightarrow \infty }nS(a/n)=S(a)$ by $\mathbb{Q}_{+}$%
-homogeneity of $S$, and further $S^{\ast }(a/n)=S(a/n)=S(a)/n\rightarrow 0,$
taking limits through $n\in \mathbb{N}.$ Taking $a<0$ gives, via Theorem 3,
that $S$ and $S^{\ast }$ are continuous. Now $S=S^{\ast },$ by continuity
and density of $\mathbb{A}$, as $S^{\ast }|\mathbb{A=}S|\mathbb{A}$. So $S|%
\mathbb{R}_{+}$ and $S|\mathbb{R}_{-}$ are linear, again by Theorem 3. The
first formulas comes from Th. HP, and the final one, in the additive case,
from $S(-t)=-S(t)$ (and then $\alpha _{S}=\beta _{S}$). $\square $

\bigskip

\noindent \textbf{Remark. }In the result above, the particular case of $S$
additive includes [BinO3, Th. 1 and Th. 1$_{\text{b}}^{\prime }].$

\bigskip

\textbf{\S 4. Complements}

\noindent 4.1 \textit{Approximate homomorphisms. }There are results in which
one has a property, such as additivity, which holds only \textit{%
approximately, }and deduces that, under suitable restrictions, it holds 
\textit{exactly. }For example, in Badora's almost-everywhere version of the
Hahn-Banach theorem [Bad4], if the relevant differences are \textit{bounded, 
}as in [Bad2]\textit{, }then they \textit{vanish. }That is, the relevant
differences are either identically plus infinity or identically zero. This
is a \textit{dichotomy, }reminiscent of those that occur in probability
theory in connection with 0-1 laws (for example, \textit{Belyaev's dichotomy 
}[Bel]; [MarR, 5.3]).

\noindent 4.2 \textit{Popa (circle) group subadditivity. }We recall from
[BinO3] that the Popa circle operation on $\mathbb{R}$, introduced in [Pop]
(cf. [Jav]), given by%
\[
a\circ b=a+b\eta (a),\text{ for }\eta (t):=1+\rho t\text{ with }\rho \geq 0, 
\]%
turns $\mathbb{G}_{+}:=\{x\in \mathbb{R}:1+\rho x>0\}$ into a group with $%
\mathbb{R}_{+}$ as a subsemigroup. The latter induces an order on $\mathbb{G}%
_{+}$ which agrees with the usual order (cf. e.g. [FucL]). So a function $f:(%
\mathbb{R}_{+},\times )\rightarrow (\mathbb{G}_{+},\circ )$ satisfying%
\[
f(xy)\leq f(x)\circ f(y) 
\]%
may be viewed as subadditive in the group context. This abstract viewpoint
encompasses both the current context of subadditivity (for $\rho =0),$ and a
further significant one arising in the theory of regular variation (the
`Goldie Functional Inequality', for $\rho =1$ -- cf. [Jab2]); for the latter
see [BinO3]. We hope to return to these matters elsewhere.

\noindent 4.3 \textit{Restricted domain. }There are results when, as in \S 3
on quantifier weakening, a property such as additivity or subadditivity
holds off some exceptional set (say, almost everywhere), and the conclusion
is also similarly restricted. This goes back to work of Hyers and Ulam
[CabC], [Bad1]. See also de Bruijn [deB], Ger [Ger].

\noindent 4.4 \textit{Origin of the Heiberg-Seneta condition. }This
condition, introduced in regular variation (see [BinGT, Th. 3.2.5],
prompting its recent study in [BinO3]), as applied to a subadditive function 
$S:\mathbb{R}\rightarrow \mathbb{R}\cup \{-\infty ,+\infty \},$ took the form%
\begin{equation}
\lim \sup\nolimits_{t\downarrow 0}S(t)\leq 0.  \tag{$HS$}
\end{equation}%
For $\mathbb{A\subseteq R}$ a dense subgroup, the assumption that $S|\mathbb{%
A}$ is linear together with $(HS)$ guaranteed not only that $S$ is
finite-valued with $S(0+)=0$, but that in fact $S$ is linear, as in Th. 5,
which relates directly to [BinGT, Th. 3.2.5].

\noindent 4.5 \textit{Examples of families of locally Steinhaus-Weil sets.}

The sets listed below are typically, though not always, members of a
topology on an underlying set.

\noindent (o) $\Sigma $ a usual (Euclidean) open set in $\mathbb{R}$ (and in 
$\mathbb{R}^{n})$ -- this is the `trivial' example;

\noindent (i) $\Sigma $ density-open subset of $\mathbb{R}$ (similarly in $%
\mathbb{R}^{n})$ (by Steinhaus's Theorem -- see e.g. [BinGT, Th. 1.1.1],
[BinO5], [Oxt, Ch. 8]);

\noindent (ii) $\Sigma $ Baire, locally non-meagre at all points $x\in
\Sigma $ (by the Piccard-Pettis Theorem -- as in [BinGT, Th. 1.1.2],
[BinO5], [Oxt, Ch. 8] -- such sets can be `thinned out', i.e. extracted as
subsets of a second-category set, using separability or by reference to the
Banach Category Theorem [Oxt, Ch.16]);

\noindent (iii) $\Sigma $ the Cantor `middle-thirds excluded' subset of $%
[0,1]$ (since $\Sigma +\Sigma =[0,2]);$

\noindent (iv) $\Sigma $ universally measurable and open in the \textit{ideal%
} topology ([LukMZ], [BinO4]) generated by omitting Haar null sets (by the
Christensen-Solecki Interior-points Theorem of [Chr1,2] and [Sol]);

\noindent (v) $\Sigma $ a Borel subset of a Polish abelian group and and
open in the ideal topology generated by omitting \textit{Haar meagre} sets
in the sense of Darji [Dar] (by Jab\l o\'{n}ska's generalization of the
Piccard Theorem, [Jab1, Th.2], cf. [Jab3], and since the Haar-meagre sets
form a $\sigma $-ideal [Dar, Th. 2.9]); for details see [BinO5].

If $\Sigma $ is \textit{Baire} (has the Baire property) and is locally
non-meagre, then it is co-meagre (since its quasi interior is everywhere
dense).

\noindent \textbf{Caveats. }1. Care is needed in identifying locally SW
sets: Mato\u{u}skov\'{a} and Zelen\'{y} [MatZ] show that in any non-locally
compact abelian Polish group there are closed \textit{non-Haar null} sets $%
A,B$ such that $A+B$ has empty interior. Recently, Jab\l o\'{n}ska [Jab4]
has shown that likewise in any \textit{non-locally compact} abelian Polish
group there are closed \textit{non-Haar meager} sets $A,B$ such that $A+B$
has empty interior.

\noindent 2. For an example on $\mathbb{R}$ of a compact subset $S$ such
that $S-S$ contains an interval, but $S+S$ has measure zero and so does not,
see [CrnGH] and the recent [BarF].

\noindent 3. Here we were concerned with subsets $\Sigma \subseteq \mathbb{R}
$ where such `anomalies' are assumed not to occur.

\noindent 4.6 \textit{Baire/measurable }$S$ and $S^{\ast }.$ Of course if $%
S\ $is Baire/measurable, then so is $S^{\ast }$, as the limsup is
sequential. Also for $\mathbb{A}$ a \textit{countable} subgroup, the upper
and lower limit functions $S_{\mathbb{A}}^{\pm }$ derived from a
subbadditive function $S$ are Baire/measurable, as the image $S(\mathbb{A})$
is countable.

\noindent 4.7 \textit{The Hahn-Banach theorem: variants. }There are various
theorems of Hahn-Banach type. Text-book accounts, as in e.g. Rudin [Rud, \S\ %
3.2, 3.3, 3.4], [FucL], deal with dominated extension theorems (without any
assumed continuity on the partial function $f$ nor on the dominating
function $p,$ HB below), separation theorems for convex sets, and continuous
extension theorems. Variations include the assumption that the dominating
function $p$ is continuous, e.g. [DodM] (implying continuity of the minorant
partial function); another variation -- from [FosM], call this `HB-lite' for
our needs in \S 4.8 below -- assumes for given $p$ merely the existence of
some linear functional dominated by $p.$ (Here, if the variant axiom is
satisfied for all $p$ continuous, then HB\ follows for all continuous $p$
[FosM, \S 4]). For a most insightful survey of very many variations in
earlier literature see [Bus]. The context also varies, correspondingly, from
vector spaces, to topological vector spaces and beyond, so to F-spaces (i.e.
topological vector spaces with topology generated by some complete
translation-invariant metric, [KalPR]) and Banach spaces. One needs to
distinguish between the variants, including the category of space over which
the assertions range, when discussing their axiomatic status. Kalton proved
([Kal], [Dre], cf. [KalPR, Ch. 4]) that an F-space in which the continuous
extension theorem (in which $f$ is continuous) holds is necessarily locally
convex, a result that is false without metrizability; it is not known
whether completeness is necessary.

The dominated extension theorem HB (i.e. without any continuity) is
equivalent to a weakened form of PIT, the Prime Ideal Theorem, namely the
existence of a non-trivial finitely-additive probability measure (as opposed
to a two-valued measure implicit in PIT) on any non-trivial Boolean algebra
([Lux], [Jec], [TomW], [MycT]) -- MB (for `measure-Boolean') in the
terminology of [MycT].

For the relative strengths of HB and the Axiom of Choice AC, see [Pin1,2];
[PinS] provide a model of set theory in which the Axiom of Dependent Choices
DC holds but HB\ fails. Moreover, HB for separable normed spaces is not
provable from DC [DodM, Cor. 4]. On the other hand, any separable normed
space satisfies the version of HB\ in which the dominating function $p$ is
continuous; indeed the partial function $f$ may first be explicitly extended
to the linear span of the union of its domain with the dense countable set
-- as in the original Banach proof [Ban] by inductive assignment of
function-values using the least possible function-value at each stage (as in
[DodM, Lemma 9]) -- and then to the rest of space, essentially as in Theorem
3, using the continuity conferred by $p$ and our sequential analysis.
(Compare [FosM] for various completeness and compactness notions here.)
Further to [DodM], we raise, and leave open here, the question as to whether
the separable case of Badura's result in [Bad3] can be proved with only DC
rather than AC, and the role that completeness (sequential or otherwise) may
play here [KalPR].

For more on axiomatics (with references), see [Bus, \S 12, 20], and Appendix
1 of the fuller arXiv version of [BinO2].

4.8. \textit{The Hahn-Banach Theorem: group analogues. }The group analogue
of the `HB-lite' property of \S 4.7 (mutatis mutandis, with `additive'
replacing `linear' etc.) delineates a class of groups providing the context
for Badora's `general' Hahn-Banach extension theorem for groups [Bad3, Th.
1], and includes amenable groups; the class is characterized in [Bad3, Th.
3] by the group analogue of HB with a side-condition on $p$. The more
special Hahn-Banach-type extension property for the case of a group $G$ of
linear operators $g:V\rightarrow V$ on a real vector space $V$ is concerned
with a $p$-dominated $G$-invariant extension of a $G$-invariant partial
linear operator $f$ (defined on a $G$-invariant subspace $W$ ) satisfying $%
f(w)\leq p(w)$ for $w\in W,$ where $p$ is a subadditive and
positive-homogeneous functional $p:V\rightarrow \mathbb{R}$ with $%
p(g(v))\leq p(v).$ This as a property of $G$ turns out to be equivalent to $%
G $ being \textit{amenable} (Silverman [Sil1,2]) -- see [Kle] for a clear
albeit early approach. See also [TomW, Th. 12.11].

\noindent 4.9 \textit{Kingman's Subadditive Ergodic Theorem. }Detailed study
of subadditivity is partially motivated by links with the Kingman
subadditive ergodic theorem, which has been very widely used in probability
theory. For background and details, see e.g. [Kin1, 2], Steele [Ste].

\bigskip

\textbf{Postscript.}

This paper germinated from the constructive and scholarly criticism of
successive drafts of [BinO3] by its Referee; it is a pleasure to thank him
again here.\newpage

\textbf{References}

\bigskip 

\noindent \lbrack Bad1] R. Badora, On some generalized invariant means and
their application to the stability of the Hyers-Ulam type. \textsl{Ann.
Polon. Math.} \textbf{LVIII} (1993), 147-159.\newline
\noindent \lbrack Bad2] R. Badora, On approximate ring homomorphisms. 
\textsl{J. Math. Anal. Appl.} \textbf{276} (2002), 589--597. \newline
\noindent \lbrack Bad3] R. Badora, On the Hahn-Banach theorem for groups. 
\textsl{Arch. Math.} (Basel) \textbf{86} (2006), 517--528.\newline
\noindent \lbrack Bad4] R. Badora, The Hahn-Banach theorem almost
everywhere. \textsl{Aequationes Math.} \textbf{90} (2016), no. 1, 173--179.%
\newline
\noindent \lbrack Ban] S. Banach, \textsl{Th\'{e}orie des op\'{e}rations lin%
\'{e}aires.} Monografie Mat. \textbf{1}, Warszawa, 1932; Chelsea Publishing
Co., New York, 1955; in \textsl{Oeuvres}, \textbf{II}, PWN-1979 (see
http://matwbn-old.icm.edu.pl/), \'{E}ditions Jacques Gabay, Sceaux, 1993; 
\textsl{Theory of linear operations.} Translated from the French by F.
Jellett. North-Holland Mathematical Library\textbf{\ 38}, 1987.\newline
\noindent \lbrack BarF] A. Bartoszewicz, M. Filipczak, Remarks on sets with
small differences and large sums. \textsl{J. Math. Anal. Appl.}\textbf{\ 456}
(2017), 245--250.\newline
\noindent \lbrack BarFN] A. Bartoszewicz, M. Filipczak, T. Natkaniec, On
Smital properties. \textsl{Topology Appl.} \textbf{158} (2011), 2066--2075.%
\newline
\noindent \lbrack Bel] Yu. K. Belyaev, Continuity and H\"{o}lder's
conditions for sample functions of stationary Gaussian processes. \textsl{%
Proc. 4th Berkeley Sympos. Math. Statist. and Prob. }\textbf{II}, 23--33
Univ. California Press, 1961.\newline
\noindent \lbrack Ber] E. Berz, Sublinear functions on $\mathbb{R},$ \textrm{%
\textsl{Aequat. Math.}} \textbf{12} (1975), 200-206.\newline
\noindent \lbrack BinGT] N. H. Bingham, C. M. Goldie and J. L. Teugels, 
\textsl{Regular variation}, 2nd ed., Cambridge University Press, 1989 (1st
ed. 1987).\newline
\noindent \lbrack BinO1] N. H. Bingham and A. J. Ostaszewski, Kingman,
category and combinatorics. \textsl{Probability and Mathematical Genetics}
(Sir John Kingman Festschrift, ed. N. H. Bingham and C. M. Goldie), 135-168,
London Math. Soc. Lecture Notes in Mathematics \textbf{378}, CUP, 2010. 
\newline
\noindent \lbrack BinO2] N. H. Bingham and A. J. Ostaszewski,
Category-measure duality: convexity, mid-point convexity and Berz
sublinearity, \textsl{Aequationes Math.}, \textbf{91.5} (2017), 801--836, (
fuller version: arXiv1607.05750).\newline
\noindent \lbrack BinO3] N. H. Bingham and A. J. Ostaszewski, Additivity,
subadditivity and linearity: automatic continuity and quantifier weakening, 
\textsl{Indagationes Mathematicae}, to appear (arXiv 1405. 3948v3).\newline
\noindent \lbrack BinO4] N. H. Bingham and A. J. Ostaszewski, Beyond
Lebesgue and Baire IV: Density topologies and a converse Steinhaus-Weil
theorem, \textsl{Topology and its Applications}, to appear; extended
version: arXiv1607.00031v2.\newline
\noindent \lbrack BinO5] N. H. Bingham and A. J. Ostaszewski, The
Steinhaus-Weil property: its converse, Solecki amenability and
subcontinuity, arXiv1607.00049v3.\newline
\noindent \lbrack Bus] G. Buskes, The Hahn-Banach theorem surveyed. \textsl{%
Dissertationes Math.} (Rozprawy Mat.) \textbf{327} (1993), 49 pp.\newline
\noindent \lbrack BusR] G. Buskes, A. van Rooij, Hahn-Banach for Riesz
homomorphisms. \textsl{Indag. Math. }\textbf{51} (1989), no. 1, 25--34.%
\newline
\noindent \lbrack CabC] F Cabello Sanchez and J. M. F. Castillo, Banach
space techniques underpinning a theory for nearly additive mappings, \textsl{%
Dissertationes Math.} \textbf{404} (2002), 73 pp.\newline
\noindent \lbrack Chr1] J. P. R. Christensen, On sets of Haar measure zero
in abelian Polish groups. Proceedings of the International Symposium on
Partial Differential Equations and the Geometry of Normed Linear Spaces
(Jerusalem, 1972). \textsl{Israel J. Math.} \textbf{13} (1973), 255--260.%
\newline
\noindent \lbrack Chr2] J. P. R. Christensen, \textsl{Topology and Borel
structure. Descriptive topology and set theory with applications to
functional analysis and measure theory.} North-Holland Mathematics Studies,
Vol. 10, 1974.\newline
\noindent \lbrack CrnGH] M. Crnjac, B. Gulja\v{s}, H. I. Miller, On some
questions of Ger, Grubb and Kraljevi\'{c}. \textsl{Acta Math. Hungar.} 
\textbf{57} (1991), 253--257.\newline
\noindent \lbrack Dar] U. B. Darji, On Haar meager sets. \textsl{Topology
Appl.} \textbf{160} (2013), 2396--2400.\newline
\noindent \lbrack deB] N. G. de Bruijn, On almost additive functions. 
\textsl{Colloq. Math.} \textbf{XV} (1966), 59-63. \newline
\noindent \lbrack DodM] J. Dodu, M. Morillon, The Hahn-Banach property and
the axiom of choice. \textsl{Math. Log. Q. }\textbf{45.3} (1999), 299--314.%
\newline
\noindent \lbrack Dre] L. Drewnowski, The weak basis theorem fails in
non-locally convex F-spaces. \textsl{Canad. J. Math.} \textbf{29} (1977),
1069--1071.\newline
\noindent \lbrack FosM] J. Fossy, and M. Morillon, The Baire category
property and some notions of compactness. J\textsl{. London Math. Soc. }(2) 
\textbf{57} (1998), 1--19.\newline
\noindent \lbrack Fuc] B. Fuchssteiner, On exposed semigroup homomorphisms. 
\textsl{Semigroup Forum} \textbf{13} (1976/77), 189--204.\newline
\noindent \lbrack FucL] B. Fuchssteiner, W. Lusky, \textsl{Convex cones. }%
North-Holland Math. Studies 56, 1981.\newline
\noindent \lbrack GajK] Z. Gajda, Z. Kominek, On separation theorems for
subadditive and superadditive functionals. \textsl{Studia Math.} \textbf{100}
(1991), 25--38.\newline
\noindent \lbrack Ger] R. Ger, On some functional equations with a
restricted domain, I, II. \textsl{Fund. Math.} \textbf{89 }(1975), 131-149; 
\textbf{98} (1978), 249-272. \newline
\noindent \lbrack HilP] E. Hille and R. S. Phillips, \textsl{Functional
analysis and semi-groups}, Coll. Publ. Vol 31, Amer. Math. Soc, 3$^{\text{rd}%
}$ ed. 1974 (1$^{\text{st}}$ ed. 1957).\newline
\noindent \lbrack Jab1] E. Jab\l o\'{n}ska, Some analogies between Haar
meager sets and Haar null sets in abelian Polish groups. \textsl{J. Math.
Anal. Appl.} \textbf{421} (2015),1479--1486.\newline
\noindent \lbrack Jab2] E. Jab\l o\'{n}ska, On solutions of a composite type
functional inequality,\textsl{\textsl{\ Math. Ineq. Apps. }}\textbf{18}
(2015), 207-215.\textsl{\textsl{\newline
}}\noindent \lbrack Jab3] E. Jab\l o\'{n}ska, A theorem of Piccard's type
and its applications to polynomial functions and convex functions of higher
orders. \textsl{Topology Appl.} \textbf{209} (2016), 46--55.\newline
\noindent \lbrack Jab4] E. Jab\l o\'{n}ska, A theorem of Piccard's type in
abelian Polish groups.\textsl{\ Anal. Math.} \textbf{42} (2016), 159--164.%
\newline
\noindent \lbrack Jav] P. Javor, On the general solution of the functional
equation $f(x+yf(x))=f(x)f(y).$ \textsl{Aequat. Math.} \textbf{1} (1968),
235-238.\newline
\noindent \lbrack Jec] T. Jech,\textsl{\ The axiom of choice.} Studies in
Logic and the Foundations of Math. \textbf{75}. North-Holland, 1973\newline
\noindent \lbrack Kal] N. J. Kalton, Basic sequences in F-spaces and their
applications. \textsl{Proc. Edinburgh Math. Soc.} (2) \textbf{19} (1974/75),
no. 2, 151--167.\newline
\noindent \lbrack KalPR] N. J. Kalton, N. T. Peck, J. W. Roberts, \textsl{An
F-space sampler.} London Math. Soc. Lect. Note Ser. \textbf{89}, Cambridge
University Press, 1984.\newline
\noindent \lbrack Kau] R. Kaufman, Extension of functionals and inequalities
on an abelian semi-group. \textsl{Proc. Amer. Math. Soc.} \textbf{17}
(1966), 83--85.\newline
\noindent \lbrack Kin1] J. F. C. Kingman, Subadditive ergodic theory. 
\textsl{Ann. Probability} \textbf{1} (1973), 883--909.\newline
\noindent \lbrack Kin2] J. F. C. Kingman, Subadditive processes. \textsl{%
Ecole d'\'{e}t\'{e} de Probabilit\'{e}s de Saint Fleur} \textsl{V}, \textsl{%
Lecture Notes in Math.} \textbf{539} (1976), 167-223.\newline
\noindent \lbrack Kle] V. Klee, Invariant extension of linear functionals. 
\textsl{Pacific J. Math.} \textbf{4} (1954), 37--46.\newline
\noindent \lbrack Kol] A. N. Kolmogoroff, Zur Normierbarkeit eines
allgemeinen topologischen linearen Raumes, \textsl{Studia Math.} \textbf{5}
(1934), 29-33.\newline
\noindent \lbrack Kra] P. Kranz, Additive functionals on abelian semigroups.%
\textsl{\ Comment. Math. Prace Mat.} \textbf{16} (1972), 239--246.\newline
\noindent \lbrack Kuc1] M. Kuczma, Almost convex functions, \textsl{Coll.
Math.} \textbf{XXI} (1070), 279-284.\newline
\noindent \lbrack Kuc2] M. Kuczma, \textsl{An introduction to the theory of
functional equations and inequalities. Cauchy's equation and Jensen's
inequality.} 2$^{\text{nd}}$ ed., Birkh\"{a}user, 2009 [1st ed. PWN,
Warszawa, 1985].\newline
\noindent \lbrack LukMZ] J. Luke\v{s}, J. Mal\'{y}, L. Zaj\'{\i}\v{c}ek, 
\textsl{Fine topology methods in real analysis and potential theory},
Lecture Notes in Mathematics \textbf{1189}, Springer, 1986.\newline
\noindent \lbrack Lux] W. A. J. Luxemburg, Reduced powers of the real number
system and equivalents of the Hahn-Banach extension theorem. in: \textsl{%
Applications of Model Theory to Algebra, Analysis, and Probability}
(Internat. Sympos., Pasadena, Calif., 1967), 123--137 Holt, Rinehart and
Winston, 1969.\newline
\noindent \lbrack MarR] M. B. Marcus and J. Rosen, \textsl{Gaussian
processes, Markov processes and local times}. Cambridge University Press,
2006.\newline
\noindent \lbrack MatZ] E. Mato\u{u}skov\'{a}, M. Zelen\'{y}, A note on
intersections of non--Haar null sets, \textsl{Colloq. Math.} \textbf{96}
(2003), 1-4.\newline
\noindent \lbrack MycT] J. Mycielski and G. Tomkowicz, Shadows of the axiom
of choice in the universe $L(\mathbb{R}).$ \textsl{Arch. Math. Logic} (doi
10.1007/s00153-017-0596-x).\newline
\noindent \lbrack Oxt] J. C. Oxtoby, \textsl{Measure and category}, 2$^{%
\text{nd}}$ ed. Grad. Texts in Math. \textbf{2}, Springer, 1980.\newline
\noindent \lbrack Pin1] D. Pincus, Independence of the prime ideal theorem
from the Hahn Banach theorem. \textsl{Bull. Amer. Math. Soc.} \textbf{78}
(1972), 766--770.\newline
\noindent \lbrack Pin2] D. Pincus, The strength of the Hahn-Banach theorem.%
\textsl{\ Victoria Symposium on Nonstandard Analysis }(Univ. Victoria,
Victoria, B.C., 1972), pp. 203--248, Lecture Notes in Math. \textbf{369},
Springer 1974.\newline
\noindent \lbrack PinS] D. Pincus, R. Solovay, Definability of measures and
ultrafilters.\textsl{\ J. Symbolic Logic} \textbf{42.2} (1977), 179--190.%
\newline
\noindent \lbrack Pop] C. G. Popa, Sur l'\'{e}quation fonctionelle $%
f[x+yf(x)]=f(x)f(y),$ \textsl{Ann. Polon. Math.} \textbf{17} (1965), 193-198.%
\newline
\noindent \lbrack Rud] W. Rudin, \textsl{Functional analysis}. 2$^{\text{ed}}
$, McGraw-Hill, 1991.\newline
\noindent \lbrack Sil1] R. J. Silverman, Invariant linear functions. \textsl{%
Trans. Amer. Math. Soc.} \textbf{81} (1956), 411--424.\newline
\noindent \lbrack Sil2] R. J. Silverman, Means on semigroups and the
Hahn-Banach extension property. \textsl{Trans. Amer. Math. Soc.} \textbf{83}
(1956), 222--237.\newline
\noindent \lbrack Sol] S. Solecki, Amenability, free subgroups, and Haar
null sets in non-locally compact groups. \textsl{Proc. London Math. Soc.} 
\textbf{(3) 93} (2006), 693--722.\newline
\noindent \lbrack Ste] J. M. Steele, Kingman's subadditive ergodic theorem, 
\textsl{Ann. de l'Institut Henri Poincar\'{e}}, \textbf{25} (1989), 93-98.%
\newline
\noindent \lbrack TomW] G. Tomkowicz and S. Wagon, \textsl{The Banach-Tarski
paradox.} Cambridge University Press, 2016 (1$^{\text{st}}$ ed. 1985).

\bigskip

\bigskip

\noindent Mathematics Department, Imperial College, London SW7 2AZ;
n.bingham@ic.ac.uk \newline
Mathematics Department, London School of Economics, Houghton Street, London
WC2A 2AE; A.J.Ostaszewski@lse.ac.uk\newpage

\end{document}